\documentclass[oneside]{amsart}
\usepackage{amsmath, amsthm,amssymb,enumerate}
\usepackage[colorlinks=true,linkcolor=black,anchorcolor=black,citecolor=black,filecolor=black,menucolor=black,runcolor=black,urlcolor=black]{hyperref}
\usepackage[T1]{fontenc}
\newcommand{\Oo}{\mathcal{O}}
\newcommand{\Z}{\mathbb{Z}}
\newcommand{\Q}{\mathbb{Q}}
\newcommand{\F}{\mathbb{F}}
\newcommand{\Tz}{\mathbb{T}_{\Z_p,\mathfrak{m}_{\Z_p}}}
\newcommand{\Jz}{J_{\Z_p,\mathfrak{m}_{\Z_p}}}
\newcommand{\To}{\mathbb{T}_{\Oo,\mathfrak{m}}}
\newcommand{\Jo}{J_{\Oo,\mathfrak{m}}}
\newcommand{\length}{\mathrm{length}}
\newcommand{\val}{\mathrm{val}}

\theoremstyle{plain}
\newtheorem{thm}{Theorem}[section]

\newtheorem{theorem}[thm]{Theorem}
\newtheorem{lemma}[thm]{Lemma}
\newtheorem{corollary}[thm]{Corollary}
\newtheorem{proposition}[thm]{Proposition}

\theoremstyle{definition}
\newtheorem{remark}[thm]{Remark}
\newtheorem{definition}[thm]{Definition}

\title{Higher congruences between newforms and Eisenstein series of squarefree level}

\author{C. Hsu}
\begin{document}
\pagestyle{plain}

\begin{abstract} Let $p\geq 5$ be prime. For elliptic modular forms of weight 2 and level $\Gamma_0(N)$ where $N>6$ is squarefree, we bound the depth of Eisenstein congruences modulo $p$ (from below) by a generalized Bernoulli number with correction factors and show how this depth detects the local non-principality of the Eisenstein ideal. We then use admissibility results of Ribet and Yoo to give an infinite class of examples where the Eisenstein ideal is not locally principal. Lastly, we illustrate these results with explicit computations and give an interesting commutative algebra application related to Hilbert--Samuel multiplicities.\end{abstract}

\maketitle

\section{Introduction}

\noindent Let $f_1,\dots,f_r$ be all weight 2 normalized cuspidal simultaneous eigenforms of level $\Gamma_0(N)$ with $N$ prime. A celebrated result of Mazur \cite[Proposition II.5.12, Proposition II.9.6]{mazur1977modular} states that if a prime $p$ divides the numerator $\mathcal{N}$ of $\frac{N-1}{12},$ then at least one of these forms is congruent modulo $p$ to the weight 2 normalized Eisenstein series \[E_{2,N}=\frac{N-1}{24}+\sum_{n=1}^\infty\sigma^\ast(n)q^n,\] where $\sigma^\ast(n)$ is the sum of all non-zero divisors $d$ of $n$ such that $(d,N)=1$. Berger, Klosin, and Kramer \cite[Proposition 3.1]{berger2014higher} refine this result to give a precise relation between $\mathrm{val}_p(\mathcal{N})$ and the depth of congruence between the newforms $f_1,\dots,f_r$ and $E_{2,N}$. Using a commutative algebra result (restated as Theorem \ref{bkkcommalg} of this paper), they show that if 
$\varpi_N$ is a uniformizer in the valuation ring of a finite extension of $\Q_p$ (of ramification index $e_N$) that contains all Hecke eigenvalues of the $f_i$'s, and $m_i$ is the largest integer such that the Hecke eigenvalues of $f_i$ and $E_{2,N}$ satisfy \[\lambda_\ell(f_i)\equiv \lambda_\ell(E_{2,N})\pmod{\varpi_N^{m_i}},\]
for all Hecke operators $T_\ell$ with $\ell\nmid N$ prime, then \begin{equation}\label{klosinineq}\frac{1}{e_N}\sum_{i=1}^r m_i\geq\mathrm{val}_p(\mathcal{N}).\end{equation}
Moreover, Theorem \ref{bkkcommalg} implies that this expression is an equality if and only if the Eisenstein ideal is locally principal.  Since the Eisenstein ideal is locally principal when $N$ is prime \cite[Theorem II.18.10]{mazur1977modular}, Eq. (\ref{klosinineq}) is always an equality in this case. However, the approach of comparing the depth of Eisenstein congruences modulo $p$, i.e., the left side of Eq. (\ref{klosinineq}), to a certain $p$-adic value suggests a way to determine if the Eisenstein ideal is locally principal for a fixed squarefree level $N$.

Let $N = \prod_{j=1}^tq_j>6$ be a squarefree positive integer. The weight 2 Eisenstein subspace of level $\Gamma_0(N),$ denoted $E_2(\Gamma_0(N)),$ is spanned by $2^t-1$ Eisenstein series, each of which is a simultaneous eigenform for all Hecke operators. Since a basis of such eigenforms can be obtained using level raising techniques \cite[\S 2.2]{yoo2013index}, each eigenform in $E_{2}(\Gamma_0(N))$ has Hecke eigenvalue $\lambda_\ell=1+\ell$ for  Hecke operators $T_\ell$ with $\ell\nmid N$ prime. Since we are interested in congruences away from $N,$ i.e., congruences between the $\ell^{th}$ Hecke eigenvalues for primes $\ell\nmid N,$ any normalized Eisenstein series $E\in E_2(\Gamma_0(N))$ will work for our generalization to squarefree level. So, let $f_1,\dots,f_r$ be all weight 2 newforms of level $\Gamma_0(N_{f_i})$ where $N_{f_i}$ divides $N$. We consider congruences between the set of newforms $f_1,\dots,f_r$ and the weight 2 Eisenstein series of level $N,$ 
\[E_{2,N}(z) = \sum_{d|N}\mu(d)dE_2(dz),\] where $\mu$ is the M\"{o}bius function and $E_2$ is the weight 2 Eisenstein series for $\mathrm{SL}(2,\Z),$ normalized so that the Fourier coefficient of $q$ is 1. Note that this Eisenstein series coincides with the original $E_{2,N}$ defined for prime levels and has $q$-expansion 
\[E_{2,N}=(-1)^{t+1}\frac{\varphi(N)}{24}+\sum_{n=1}^\infty\sigma^\ast(n)q^n.\]

Now, in contrast Mazur's result for the prime level case, Ribet--Yoo \cite{yoo2014non} establish that Eisenstein congruences can occur for primes $p$ such that some prime divisor $q_j$ of $N$ satisfies $q_j\equiv\pm 1\pmod p$. We therefore consider a function $\eta(N)$ given by \[\eta(N) =\prod_{j=1}^t(q_j^2-1) = \varphi(N)\cdot\prod_{j=1}^t(q_j+1).\]

The first main result of this paper (Proposition \ref{congruenceresult}) extends the higher congruences framework of Berger, Klosin, and Kramer to squarefree level $N>6$ so that after replacing $\val_p(\mathcal{N})$ with $\val_p(\eta(N)),$ the inequality in Eq. (\ref{klosinineq}) still holds. Under some mild assumptions on $N$ and $p,$ the second main result (Theorem \ref{princondition}) then gives a numerical criterion, in terms of the depth of Eisenstein congruences, for the Eisenstein ideal to not be locally principal. The last main result, stated below, uses this numerical criterion and the existence of sufficiently many Eisenstein congruences to give a condition,  in terms of only $\val_p(\eta(N)),$ for the Eisenstein ideal to not be locally principal: 
\begin{theorem}\label{nonprincthm} Let $N =\prod_{i=1}^t q_i$ be a squarefree integer and $p\geq 5$ be a prime such that $p\nmid N$. Assume $\val_p(q_i-1)=0$ for all $i,$ and $\val_p(q_i+1)>0$ for $i = 1,\dots,s,$ where $1\leq s\leq t$. If $s\cdot 2^{t-2}>\val_p(\eta(N)),$ then $J_\Z$ is not locally principal.
\end{theorem}

As an application of these results, we express the depth of congruence \[\frac{1}{e_N}\sum_{i=1}^rm_i\] (from Proposition \ref{congruenceresult}) as the Hilbert--Samuel multiplicity of the Eisenstein ideal in the Hecke algebra. While the depth of Eisenstein congruences modulo $p$ detects whether the associated (local) Eisenstein ideal is principal, this connection to multiplicities might allow us to give a more precise statement regarding the minimal number of generators.

Using an algorithm adapted from Naskr\k{e}cki \cite[\S 4.2]{naskrkecki2014higher}, we provide computational examples to illustrate our main results. While Naskr\k{e}cki has computed a large number of Eisenstein congruences, his work concerns congruences of $q$-expansions rather than congruences away from $N$. As a result, his data does not necessarily agree with ours. For example, if $N=97$ and $p=2,$ then $\mathrm{val}_2(\frac{96}{12}) = 3$. Since the constant term of $E_{2,N}$ has a $2$-adic valuation of 2, Naskr\k{e}cki's algorithm returns 2 as the depth of congruence. On the other hand, our algorithm returns 3, which agrees with the equality in Eq. (\ref{klosinineq}).  Moreover, Naskr\k{e}cki's algorithm determines the exact Eisenstein series in $E_2(\Gamma_0(N))$ for which a congruence holds; we do not require this information since the Hecke eigenvalues of all Eisenstein series in $E_2(\Gamma_0(N))$ coincide away from $N$. Because of these differences, we use our modified algorithm for congruence computations. 

Lastly, we note that work of Wake and Wang-Erickson \cite{wake2017rank,wake2018eisenstein}, Ribet--Yoo \cite{yoo2013index}, and Ohta \cite{masami2014eisenstein} independently study similar questions about various Eisenstein ideals. 

This paper is organized as follows. In Section \ref{commalg}, we recall important commutative algebra results as well as background information on the Hilbert--Samuel function. In Section \ref{depthsresult}, we prove the main results and briefly address the cases $p=2,3$. In Section \ref{examples and application}, we discuss applications and give computational examples. Appendix \ref{algorithm} contains the algorithm used to compute all congruences.

\subsection*{Acknowledgments} This paper grew out of a PhD thesis completed at the University of Oregon under the supervision of Ellen Eischen, whom I wish to thank for her constructive feedback and guidance throughout. I am also grateful to Krzysztof Klosin for generously taking time to provide helpful insights and detailed comments on this project during the past year. I thank Tobias Berger, Kimball Martin, Preston Wake, and Carl Wang-Erickson for several useful conversations. Lastly, I thank the referee for a careful review of the manuscript. This work was partially supported by an AAUW American Dissertation Fellowship and a UO Doctoral Research Fellowship.

\section{Commutative algebra preliminaries}\label{commalg}
Throughout this section, we use the following notation. Let $p$ be a prime, and let $\mathcal{O}$ be the valuation ring of a finite extension $E$ of $\Q_p$. Also, let $\varpi$ be a uniformizer of $\Oo$ and write $\F_\varpi=\Oo/\varpi\Oo$ for the residue field. 

For $s\in\Z_+,$ let $\{n_1,\,n_2,\,\dots,\,n_s\}$ be a set of $s$ positive integers and set $n = \sum_{i=1}^sn_i$. For each $i\in\{1,\,2,\,\dots,\,s\},$ let $A_i = \Oo^{n_i}$ and set $A=\prod_{i=1}^sA_i=\Oo^n$. Also, let $T\subset A$ be a local complete $\Oo$-subalgebra which is of full rank as an $\Oo$-submodule, and let $J\subset T$ be an ideal of finite index. For each $i,$ we define $\varphi_i:A\twoheadrightarrow A_i$ to be the canonical projection and set $T_i=\varphi_i(T)$ and $J_i=\varphi_i(J)$. Note that since each $T_i$ is also a (local complete) $\Oo$-subalgebra and the projections $\varphi_i|_T$ are local homomorphisms, $J_i$ is also an ideal of finite index in $T_i$. 

We first recall a result of Berger, Klosin, and Kramer \cite[Theorem 2.1]{berger2014higher} which is key to proving Proposition \ref{congruenceresult}. We then define the Hilbert--Samuel function of the module $T$ as well as the associated multiplicity $e(J,T)$ of the ideal $J\subset T$. In particular, we prove that \[e(J,T) =\sum_{i=1}^s \mathrm{length}(T_i/J_i).\]

\subsection{Result of Berger, Klosin, and Kramer} Using the Fitting ideal $\mathrm{Fit}_\Oo(M)$ associated to a finitely presented $\Oo$-module $M$ (cf. \cite[Appendix]{mazur1984class}), Berger, Klosin, and Kramer prove the following commutative algebra result, which is widely applicable in the context of congruences between automorphic forms:

\begin{theorem}[Berger--Klosin--Kramer, 2013]\label{bkkcommalg} If $\#\mathbb{F}_\varpi^\times\geq s-1$ and each $J_i$ is principal, then \[\#\prod_{i=1}^sT_i/J_i\geq\#T/J.\] Moreover,  the ideal $J$ is principal if and only if equality holds.
\end{theorem}

\begin{remark} Note that this inequality is often strict, as illustrated in the application of Theorem \ref{bkkcommalg} to Eisenstein congruences of elliptic modular forms of squarefree level.
\end{remark}

\subsection{The Hilbert--Samuel function and multiplicities} 
Let $R$ be a local ring with maximal ideal $\mathfrak{m}$. For a finitely generated $R$-module $M$ and an ideal $\mathfrak{q}\subset R$ of finite colength on $M,$ define the \textit{Hilbert--Samuel function} of $M$ with respect to $\mathfrak{q}$ to be (cf. \cite[\S 12.1]{eisenbud2013commutative}) \[H_{\mathfrak{q},M}(n):=\mathrm{length}(\mathfrak{q}^nM/\mathfrak{q}^{n+1}M).\]
By \cite[Theorem 12.4]{eisenbud2013commutative}, we have \[\dim M = 1+\deg P_{\mathfrak{q},M},\] where $P_{\mathfrak{q},M}(n)$ is a polynomial that agrees with $H_{\mathfrak{q},M}(n)$ for large enough $n$. 

Moreover, by \cite[Exercise 12.6]{eisenbud2013commutative}, we may write
\[P_{\mathfrak{q},M}(n) = \sum_{i=0}^da_iF_i(n),\]
where $F_i(n) = \left(\substack{n\\i}\right)$ is the binomial coefficient regarded as a polynomial in $n$ of degree $i,$ and the $a_i$ are integers with $a_d\neq 0$. Given these functions, we have the following definition:

\begin{definition}
The coefficient $a_d$ is called the \textit{multiplicity} of $\mathfrak{q}$ on $M$ and is denoted $e(\mathfrak{q},M)$.
\end{definition} 

Note that the leading coefficient of $P_{\mathfrak{q},M}$ equals $e(\mathfrak{q},M)/d!$. In particular, when $M$ is a finitely generated free $R$-module and $\dim R=1,$ we have $\dim M:=\dim R/\mathrm{Ann}_R{M} = 1,$ and hence, $P_{\mathfrak{q},M
}$ will be a constant function for large enough $n$. Thus, in this case,  \[P_{\mathfrak{q},M} = e(\mathfrak{q},M).\]

To relate the multiplicity $e(J,T)$ to $\sum_i\length(T_i/J_i),$ we use the following proposition:

\begin{proposition}\label{multiprop}If each $J_i$ is principal, then we have
\begin{equation}\label{mult}e(J,T) = \sum_{i=1}^s e(J_i,T_i).\end{equation}
\end{proposition}

\begin{remark}\label{elength} When $J=(\alpha)$ is principal, this equality follows immediately from \cite[Proposition 2.3]{berger2014higher}. Indeed, since multiplication by $\alpha$ gives $T$-module isomorphisms 
\[T/J\simeq J/J^2\simeq J^2/J^3\simeq\cdots,\] $H_{J,T}(n)$ is a constant fuction and $e(J,T)=\mathrm{length}(T/J)$. Similarly, we have $e(J_i,T_i)=\mathrm{length}(T_i/J_i)$. 
Additionally, note that for any $J,$ 

\[\sum_{i=1}^n\mathrm{length}\left(J_i^r/J_i^{r+1}\right) = \length\left(\prod_{i=1}^n J_i^r/J_i^{r+1}\right)=\length\left(\frac{\prod_{i=1}^n J_i^r}{\prod_{i=1}^n J_i^{r+1}}\right),\] and hence,
\begin{equation}\label{directsum}\sum_{i=1}^n e(J_i,T_i) = e\left(\prod_{i=1}^n J_i,\prod_{i=1}^n T_i\right).\end{equation}
\end{remark}

We now prove Proposition \ref{multiprop} using the following two lemmas:
\begin{lemma}[Properties of Multiplicities {[5, Exercise 12.11.a.ii]}]\label{properties} Let \[0\to\,M'\to\,M\to\,M''\to\,0\]
be an exact sequence of modules over the local ring $(R,\mathfrak{m}),$ and suppose that $\mathfrak{q}\subset R$ is an ideal of finite colength on $M,\,M',\,M''$. If $\dim M=\dim M'>\dim M'',$ then $e(\mathfrak{q},M) = e(\mathfrak{q},M')$.
\end{lemma}

\begin{lemma}\label{product} We have \[J\prod_{i=1}^s T_i\subseteq \prod_{i=1}^sJ_i,\] with equality whenever the $J_i$ are principal.
\end{lemma}

\begin{proof} The left-hand side consists of elements of the form \[\alpha\cdot\left(\varphi_1(t_1),\dots,\varphi_s(t_s)\right) =\left(\varphi_1(\alpha)\varphi_1(t_1),\dots,\varphi_s(\alpha)\varphi_s(t_s)\right) \] with $\alpha\in J$ and $t_j\in T,$ and hence, the containment is clear since $\varphi_i(J) = J_i$ is an ideal of $T_i$. When the $J_i$ are principal, \cite[Proposition 2.6]{berger2014higher} guarantees the existence of some $\alpha\in J$ such that $\varphi_i(\alpha)$ generates $J_i$ for all $i$. Thus, we may write an element of the right-hand side as \[\left(\varphi_1(\alpha)\varphi_1(t_1),\dots,\varphi_s(\alpha)\varphi_s(t_s)\right) = \alpha\cdot\left(\varphi_1(t_1),\dots,\varphi_s(t_s)\right)\] for some $t_1,\dots,t_s\in T$.
\end{proof}

\begin{proof}[Proof of Proposition \ref{multiprop}] Consider the exact sequence \[0\to\,T\to\,\prod_{i=1}^s T_i\to\,K\to\,0,\] where $K$ denotes the cokernel.  Since $T$ has full rank in $A,$ Lemma \ref{properties}  gives
\[\label{prodT}e(J,T)= e\left(J,\prod_{i=1}^s T_i\right),\] and hence, since the $J_i$ are principal, we can apply Lemma \ref{product} to obtain \[e\left(J,T\right)=e\left(J,\prod_{i=1}^sT_i\right)=e\left(J\prod_{i=1}^s T_i,\prod_{i=1}^s T_i\right) = e\left(\prod_{i=1}^{s}J_i,\prod_{i=1}^sT_i\right).\]
Thus, by Eq. (\ref{directsum}), \[e(J,T) = \sum_{i=1}^s e(J_i,T_i).\] 
\end{proof}

\begin{corollary}\label{length} If each $J_i$ is principal, then \[e(J,T)=\sum_{i=1}^s\length(T_i/J_i).\]\end{corollary}

\begin{proof} As established in Remark \ref{elength}, if each $J_i$ is principal, $e(J_i,T_i) = \mathrm{length}(T_i/J_i)$. Hence, \[e(J,T) =\sum_{i=1}^s e(J_i,T_i)=\sum_{i=1}^s \mathrm{length}(T_i/J_i).\]
\end{proof}

\section{Higher congruences: Proof of main results}\label{depthsresult}
Recall that $N = \prod_{j=1}^tq_j>6$ is a squarefree positive integer,  $f_1,\dots,f_r$ are all weight 2 newforms of level $N_{f_i}$ dividing $N,$ and  \[E_{2,N}=(-1)^{t+1}\frac{\varphi(N)}{24}+\sum_{n=1}^\infty\sigma^\ast(n)q^n.\] In this section, we first extend the higher congruences framework in \cite[\S 3]{berger2014higher} to elliptic modular forms of squarefree level. Under certain conditions, we then give two numerical criteria, one in terms of the depth of Eisenstein congruences modulo $p$ and one in terms of the $p$-valuation of $\eta(N),$ for the Eisenstein ideal to not be locally principal. In particular, these criteria allow us to establish an infinite class of examples where the Eisenstein ideal is not locally principal.
\subsection{Higher congruences framework for squarefree level}\label{framekwork} For each prime $p\geq 5,$ we would like to bound the depth of Eisenstein congruences modulo $p$ by the $p$-adic valuation of the index of an Eisenstein ideal in the associated Hecke algebra. Indeed, there are many choices for which Hecke algebra to study. In our context, we are only interested in congruences away from $N$, and so we consider Hecke algebras generated by the Hecke operators $T_\ell$ for primes $\ell\nmid N$. We refer to such Hecke algebras as anemic in order to  emphasize the exclusion of the Hecke operators $T_q$ for primes $q\,|\,N$. (When the level $N$ is prime, it makes no difference whether we include $T_N$ since this Hecke operator acts as the identity in the associated Hecke algebra \cite[Proposition 3.19]{calegari2005ramification}.) Additionally, because we want to use arithmetic data from certain Galois representations, it can be convenient to exclude the Hecke operator $T_p$. To distinguish between whether we exclude or include $T_p,$ we write $\mathbb{T}$ to denote the Hecke algebra generated by $T_\ell$ for primes $\ell\nmid Np$ and $\tilde{\mathbb{T}}$ to denote the Hecke algebra generated by $T_\ell$ for primes $\ell\nmid N$. Note that while we can formulate most of the results in this section for either $\mathbb{T}$ or $\tilde{\mathbb{T}},$ for simplicity's sake, we mainly use the Hecke algebra $\mathbb{T}$.

Let $S_2(N)$ denote the $\mathbb{C}$-space of modular forms of weight 2 and level $\Gamma_0(N)$. For any subring $R\subset\mathbb{C},$ we write $\mathbb{T}_R$ for the $R$-subalgebra of $\mathrm{End}_\mathbb{C}(S_2(N))$ generated by the Hecke operators $T_\ell$ for primes $\ell\nmid Np$. Let $J_R$ be the Eisenstein ideal, i.e., the ideal of $\mathbb{T}_R$ generated by $T_\ell-(1+\ell)$ for primes $\ell\nmid Np$. For a prime ideal $\mathfrak{a}$ of $\mathbb{T}_R,$ write $\mathbb{T}_{R,\mathfrak{a}} = \displaystyle\varprojlim_{m}\mathbb{T}_R/\mathfrak{a}^m$ for the completion of $\mathbb{T}_R$ at $\mathfrak{a},$ and set $J_{R,\mathfrak{a}}:= J_R\mathbb{T}_{R,\mathfrak{a}}$. We will call $J_{R,\mathfrak{a}}$ the local Eisenstein ideal.

It is well-known that $S_2(N)$ is isomorphic to $\mathrm{Hom}_\mathbb{C}(\mathbb{T}_\mathbb{C}(N),\mathbb{C}),$ where $\mathbb{T}_\mathbb{C}(N)$ is the full Hecke algebra in $\mathrm{End}_\mathbb{C}(S_2(N))$ \cite[Proposition 12.4.13]{diamond1995modular}. We now establish an analogous duality between the anemic Hecke algebra $\mathbb{T}_\mathbb{C}$ and the $\mathbb{C}$-subspace $L$ of $S_2(N)$ spanned by newforms $f_1,\dots,f_r$. Consider the bilinear pairing 
\begin{equation}\label{perfectpairing}\begin{split}\mathbb{T}_\mathbb{C}\times L&\to\,\mathbb{C}\\(T,f)&\mapsto a_1(Tf),\end{split}\end{equation} where $a_n(Tf)$ denotes the $n^{th}$ Fourier coefficient of $Tf$. This pairing induces maps \begin{equation}
\label{pairingmaps}\begin{split}L&\to\,\mathrm{Hom}_\mathbb{C}(\mathbb{T}_\mathbb{C},\mathbb{C}) = \mathbb{T}_\mathbb{C}^\vee\\
\mathbb{T}_\mathbb{C}&\to\,\mathrm{Hom}_\mathbb{C}(L,\mathbb{C}) = L^\vee.
\end{split}
\end{equation}

\begin{proposition}
\label{pairingprop} The above maps are isomorphisms.
\end{proposition}

\begin{proof} Since a finite dimensional vector space and its dual have the same dimension, it suffices to show that each map is injective. To show these maps are injective, we require the following lemma, which uses Atkin--Lehner theory:
\begin{lemma}\label{pairinglemma}
Any $f\in L$ with $a_n(f) = 0$ for all $(n,Np)=1$ is $0$. 
\end{lemma}
\begin{proof} Consider $f\in L$ such that $a_n(f) = 0$ for all $(n,Np) = 1$. By \cite[Proposition 6.2.1]{diamond1995modular} or \cite[Theorem 4.6.8]{miyake2006modular}, there exist cusp forms $g_q(z)\in S_2(N/q)$ for all prime factors $q$ of $N$ such that \begin{equation}\label{eq1}f(z) = \sum_{q|N}g_q(qz).\end{equation} In particular, by Atkin--Lehner theory, cf. \cite[\S 6]{diamond1995modular}, the right-hand side of this equation can be expressed as a linear combination (over $\mathbb{C}$) of cusp forms $\{f_i(d_i^jz)\},$ where each $f_i(z)$ is a newform of level $M_i$ dividing $N,$ and for each $i,$ $d_i^j$ runs over all integers that are strictly greater than $1$ and satisfy $d_i^jM_i\,|\,N$. 

On the other hand, since $f\in L,$ we can write \begin{equation}\label{eq2}f(z)=\sum_{i=1}^rb_if_i(z),\,\,\,\,b_i\in\mathbb{C},\end{equation} and so comparing Eqs. (\ref{eq1}) and (\ref{eq2}) gives a linear dependence between the newforms $f_1(z),\dots,f_r(z)$ and the cusp forms $\{f_i(d_i^jz)\}$. However, any such dependence must be trivial since the collection \[\{f_1(z),\dots,f_r(z),f_1(d_1^1z),\dots,f_1(d_1^{j_1}z),\dots,f_r(d_r^1z),\dots,f_r(d_r^{j_r}z)\}\]forms a basis for $S_2(N),$ and so we conclude that $b_1 = \cdots=b_r = 0,$ i.e., $f= 0$.

\end{proof}

Given this lemma, we prove the injectivity of the maps in Eq. (\ref{pairingmaps}) as follows. First, suppose that $f\mapsto 0\in\mathrm{Hom}_\mathbb{C}(\mathbb{T}_\mathbb{C},\mathbb{C})$. Then $a_1(Tf) =0$ for all $T\in\mathbb{T}_\mathbb{C},$ so $a_n = a_1(T_nf) = 0$ for all $(n,Np) = 1$. By Lemma \ref{pairinglemma}, $f =0$.

Next, suppose $T\mapsto 0\in\mathrm{Hom}_\mathbb{C}(L,\mathbb{C})$ so that $a_1(Tf) = 0$ for all $f\in L$. Substituting $T_nf$ for $f$ and using the commutativity of $\mathbb{T}_\mathcal{O},$ we obtain \[a_1(T(T_nf)) =a_1(T_n(Tf)) = a_n(Tf) = 0,\] for all $(n,Np) = 1$. Hence, Lemma \ref{pairinglemma} implies that $Tf = 0$. Since any $g\in S_2(N)$ can be written as a linear combination of $\mathbb{T}_\mathbb{C}$-eigenforms \cite[Proposition 1.20]{darmon1995fermat}, and since each of these eigenforms shares its eigencharacter with some $f\in L$ \cite[Theorem 1.22]{darmon1995fermat}, we conclude that $Tg = 0$ for all $g\in S_2(N),$ i.e., $T=0$.
\end{proof}

We now apply Theorem \ref{bkkcommalg} to Eisenstein congruences of elliptic modular forms of squarefree level. Fix an embedding $\overline{\Q}_p\hookrightarrow\mathbb{C}$ and let $E$ be a finite extension of $\Q_p$ that contains all Hecke eigenvalues of the $f_i$'s and whose residue field has order at least $s$. Write $\mathcal{O}_N$ for the ring of integers in $E,$ $\varpi_N$ for a choice of uniformizer, $e_N$ for the ramification index of $\mathcal{O}_N$ over $\Z_p,$ and $d_N$ for the degree of its residue field over $\F_p$. We are interested in the local structure of the Eisenstein ideal when we complete $\mathbb{T}_{\Z_p}$ at the unique maximal ideal $\mathfrak{m}_{\Z_p}\subseteq\mathbb{T}_{\Z_p}$ containing $J_{\Z_p}$. Indeed, the following result relates the depth of Eisenstein congruences modulo $p$ to the $p$-adic valuation of $\#\mathbb{T}_{\mathbb{Z}_p}/J_{\Z_p}$ and shows that this depth detects whether the local Eisenstein ideal $\Jz$ is principal:

\begin{proposition}\label{congruenceresult}
For $i=1,\dots,r,$ let $\varpi_N^{m_i}$ be the highest power of $\varpi_N$ such that the Hecke eigenvalues of $f_i$ are congruent to those of $E_{2,N}$ modulo $\varpi_N^{m_i}$ for Hecke operators $T_\ell$ for all primes $\ell\nmid Np$. Then, we have  \begin{equation}\label{inequality}\frac{1}{e_N}(m_1+\cdots+m_r)\geq\mathrm{val}_p(\#\mathbb{T}_{\Z_p}/J_{\Z_p}).\end{equation} This inequality is an equality if and only if the Eisenstein ideal $\Jz$ is principal.
\end{proposition}

\begin{proof} To simplify notation, write $\mathcal{O}$ for $\mathcal{O}_N$ and $\varpi$ for $\varpi_N,$ and let $\mathfrak{m}=J_\Oo+\varpi\mathbb{T}_\Oo$ be the unique maximal ideal of $\mathbb{T}_\mathcal{O}$ containing $J_\mathcal{O}$. By Atkin--Lehner theory, each newform $f_1,\dots,f_r$ (of level $N_{f_i}$) is a simultaneous eigenform under the action of the anemic Hecke algebra $\mathbb{T}_\mathcal{O},$ and so we can consider the map
\begin{equation}\label{heckemap}\mathbb{T}_\Oo\to\,\prod_{i=1}^s\Oo,\,\,\,T_\ell\mapsto\,\prod_{i=1}^s(\lambda_\ell(f_i)).\end{equation}
In particular, the perfect pairing established by Proposition \ref{pairingprop} implies that this map is an injection. Indeed, if $T\in\mathbb{T}_\Oo$ maps to 0, then by viewing $T$ as a $\mathbb{C}$-linear form on $L$ via an extension of scalars, we see that $T\mapsto 0\in L^\vee$ in Eq. (\ref{pairingmaps}), i.e., $T=0$.

Now, renumber $f_1,\dots,f_r$ so that $f_1,\dots,f_s$ satisfy an Eisenstein congruence away from $N$ while $f_{s+1},\dots,f_r$ do not. Eq. (\ref{heckemap}) then induces an injection \[\To\hookrightarrow\prod_{i=1}^s\Oo,\,\,\,T_\ell\mapsto\,\prod_{i=1}^s(\lambda_\ell(f_i)).\]
Since $\To\subset\prod_{i=1}^s\Oo$ is a local complete $\Oo$-subalgebra of full rank, we apply Theorem \ref{bkkcommalg} with $T=\To,$ $J=\Jo,$\,\,$T_i = \mathcal{O},$ and $\varphi_i:T\to\,T_i$ as the canonical projection. (Note that by construction, $E$ satisfies the hypothesis in Theorem \ref{bkkcommalg} on the order of its residue field.) For each projection $T_i/J_i,$ we have \begin{equation}\label{TiJi}\val_p(\#T_i/J_i) = \val_p\left(\#\,\mathcal{O}/\varpi^{m_i}\mathcal{O}\right) = m_id_N = m_i\frac{[\Oo:\Z_p]}{e_N}.\end{equation}
On the other hand, we have $\To = \Tz\otimes_{\Z_p}\Oo$ \cite[Lemma 3.27 and Proposition 4.7]{darmon1995fermat} and $\Jo = \Jz\otimes_{\Z_p}\Oo,$ and hence,
\begin{equation}\label{scaledindex}\val_p(\#T/J)=\val_p\hspace{-.075cm}\left(\#\frac{\Tz}{\Jz}\otimes_{\Z_p}\Oo\right) =[\Oo:\Z_p]\cdot\val_p\hspace{-.075cm}\left(\#\frac{\Tz}{\Jz}\right).\end{equation}
Combining these equalities yields \[\frac{1}{e_N}(m_1+\cdots+m_r)\geq\val_p\hspace{-.075cm}\left(\#\frac{\Tz}{\Jz}\right),\] and hence, the result follows from the fact that $\Tz/\Jz\simeq\mathbb{T}_{\Z_p}/J_{\Z_p}$.
\end{proof}

\subsection{Local principality of the Eisenstein ideal for squarefree level}\label{localprincipality} To use Proposition \ref{congruenceresult} to generate examples of squarefree levels for which the Eisenstein ideal is not locally principal, we need to (i) determine the $p$-adic valuation of $\#\mathbb{T}_{\Z_p}/J_{\Z_p},$ ideally in terms of a related $L$-value, and (ii) show that the depth of Eisenstein congruence modulo $p$ is strictly greater than this $p$-adic valuation. 
\subsubsection{The index of the Eisenstein ideal inside the associated Hecke algebra} Many of the current methods used to compute the size of a congruence module attached to an Eisenstein series center on deformation theory and an $R=\mathbb{T}$ argument, c.f. \cite{berger2018modularity,berger2013deformation,skinner1997ordinary}. Specifically, one can study deformations of  mod $p$ Galois representations of dimension 2 whose semi-simplification is the direct sum of two characters. 

Since we are concerned with congruences between Eisenstein series and cusp forms of weight 2 and trivial Nebentypus, we consider mod $p$ Galois representations whose semi-simplification is the direct sum of the trivial character and the mod $p$ reduction of the $p$-adic cyclotomic character. Indeed, Berger--Klosin \cite{berger2018modularity} prove that the order of a certain Selmer group bounds the size of the congruence module. They then use the Main Conjecture of Iwasawa theory \cite{mazur1984class} to bound the order of the relevant Selmer group by a generalized Bernoulli number with correction factors, which in our setting is equal to \[\eta(N) =\prod_{j=1}^t(q_j^2-1) = \varphi(N)\cdot\prod_{j=1}^t(q_j+1).\]
Due to technical obstacles arising in their method, Berger--Klosin assume $p\nmid N$ and that each prime divisor $q_j$ of $N$ satisfies $q_j\not\equiv 1\pmod p,$ and so, for the remainder of this section, we assume that the squarefree level $N$ satisfies these conditions. 

We now state the result of Berger--Klosin that bounds $\#\mathbb{T}_{\Z_p}/J_{\Z_p}$:

\begin{proposition}[Berger--Klosin, 2018]\label{indexBK} One has \[\val_p(\eta(N))\geq\val_p(\#\mathbb{T}_{\Z_p}/J_{\Z_p}).\]
\end{proposition}

\begin{proof} This bound follows from Propositions 3.10 and 5.7 in \cite{berger2018modularity}. Note that while
 the results in \cite{berger2018modularity} concern the index $\#\To/\Jo,$ we can use
Eq. (\ref{scaledindex}) to give equivalent statements for $\#\Tz/\Jz$.
\end{proof}

Thus, Propositions \ref{congruenceresult} and \ref{indexBK} yield the following theorem:

\begin{theorem}\label{princondition} Let $N$ be a squarefree integer such that none of its prime divisors are congruent to $1\pmod p$ for a prime $p\nmid N$. If the depth of Eisenstein congruences mod $p$ is strictly greater than $\mathrm{val}_p(\eta(N)),$ i.e., \[\frac{1}{e_N}(m_1+
\cdots+m_r)>\mathrm{val}_p(\eta(N)),\] then the local Eisenstein ideal $\Jz$ is non-principal.
\end{theorem}

Now, while the reverse bound in Proposition \ref{indexBK} is not required to show that the local Eisenstein ideal $\Jz$ is non-principal, it is still of interest, particularly within the context of the modularity of residual Galois representations. Indeed, because we are considering Hecke algebras and Eisenstein ideals associated to weight 2 cusp forms with trivial Nebentypus, the methods of Berger--Klosin referenced in Proposition 
\ref{indexBK} cannot be used to establish the reverse bound (see \cite[\S 5]{berger2018modularity}). Nonetheless, under the additional assumption that at least one prime divisor of $N$ is not congruent to $-1 \mod p,$ we can use work of Ohta \cite{masami2014eisenstein} on congruence modules attached to Eisenstein series to obtain the desired lower bound. Note that Ohta includes the Hecke operator $T_p$ in his congruences modules, and so in what follows, we initially work with the Hecke algebra $\tilde{\mathbb{T}}_{\Z_p}$ (which includes the Hecke operator $T_p$) and then pass back to our usual Hecke algebra  $\mathbb{T}_{\Z_p}$ (which does not).

Following \cite[\S 2-3]{masami2014eisenstein}, we consider the action on $S_2(N)$ of the Atkin--Lehner involutions $w_d$ for all positive divisors $d$ of $N$. 
More specifically, for $N = \prod_{j=1}^t q_j,$ set $\textbf{\textit{E}} = \{\pm 1\}^t$. Then for each $\boldsymbol{\varepsilon} = (\varepsilon_1,\dots,\varepsilon_t)\in\textbf{\textit{E}},$ define $S_2(N)^{\boldsymbol{\varepsilon}}$ to be the maximum direct summand of $S_2(N)$ on which $w_{q_j}$ acts as multiplication by $\varepsilon_{j}$ $(1\leq j\leq t)$. Since the Atkin--Lehner operators $w_d$ commute with the Hecke operators $T_\ell$ for $\ell\nmid N,$ the subspace $S_2(N)^{\boldsymbol{\varepsilon}}$ is invariant under the action of $\tilde{\mathbb{T}}_{\Z_p}$. So, let $\tilde{\mathbb{T}}_{\Z_p}^{\boldsymbol{\varepsilon}}$ (resp. $\tilde{J}_{\mathbb{Z}_p}^{\boldsymbol{\varepsilon}}$) denote the restriction of $\tilde{\mathbb{T}}_{\Z_p}$ (resp. $\tilde{J}_{\Z_p}$) to $S_2(N)^{\boldsymbol{\varepsilon}}$. Here $\tilde{J}_{\Z_p}\subset\tilde{\mathbb{T}}_{\Z_p}$ denotes the Eisenstein ideal which includes the additional generator $T_p-(1+p)$.

We would like to use a result of Ohta that computes the $p$-adic valuation of the index of the Eisenstein ideal inside of a certain Hecke algebra. Since Ohta's notation differs significantly from ours, we briefly explain his notation and how it relates to our conventions. Indeed, in his work on Eisenstein ideals and rational torsion subgroups, Ohta studies three different spaces of modular forms which he denotes $M_k^\mathrm{A}(\Gamma_0(N);\Z_p),\,M_k^\mathrm{B}(\Gamma_0(N);\Z_p),$ and $M_k^\mathrm{reg}(\Gamma_0(N);\Z_p)$. The first (resp. the second) space consists of modular forms in the sense of Deligne--Rapoport and Katz (resp. Serre and Swinnerton--Dyer), and the third space consists of regular differentials on the modular curve. Since $\Z_p$ is flat over $\Z[1/N],$ these three spaces of modular forms coincide \cite[Eq. (1.3.4) and Cor 1.4.10]{masami2014eisenstein}, and so, although Ohta defines his Hecke algebra $\mathbb{T}(N;\Z_p)$ as a subring of $\mathrm{End}_{\Z_p}(S_2^{\mathrm{reg}}(\Gamma_0(N);\Z_p),$ we can view it as subring of $\mathrm{End}_{\Z_p}(S_2(\Gamma_0(N)))$.  Moreover, while $\mathbb{T}(N;\Z_p)$ includes the Atkin--Lehner involutions and so differs in general from the anemic Hecke algebra $\tilde{\mathbb{T}}_{\Z_p},$ its restriction $\mathbb{T}(N;\Z_p)^{\boldsymbol{\varepsilon}}$ to $ S_2(\Gamma_0(N))^{\boldsymbol{\varepsilon}}$ is generated by the Hecke operators $T_\ell$ for primes $\ell\nmid N$ and therefore coincides with $\tilde{\mathbb{T}}_{\Z_p}^{\boldsymbol{\varepsilon}}$. Thus, when ${\boldsymbol{\varepsilon}}\neq{\boldsymbol{\varepsilon}}_+,$ where $\boldsymbol{\varepsilon}_+ = (1,1,1,\dots,1),$ we may apply \cite[Theorem 3.1.3]{masami2014eisenstein} to obtain the equality
\begin{equation}\label{ohtaequality}\mathrm{val}_p(\#\tilde{\mathbb{T}}_{\Z_p}^{\boldsymbol{\varepsilon}}/\tilde{J}_{\mathbb{Z}_p}^{\boldsymbol{\varepsilon}}) = \mathrm{val}_p\Big(\prod_{j=1}^t(q_j+\varepsilon_j)\Big).\end{equation} In particular, since $p\neq2,$ we can choose some $\boldsymbol{\varepsilon}' = (\varepsilon_1',\dots,\varepsilon_t')\in\mathbf{E}$ such that \[\val_p(q_j+\varepsilon_j') = \val_p(q_j^2-1)\] for each $j=1,\dots,t$. Then, under the additional assumption that $q_j\not\equiv-1\pmod p$ for at least one $j,$ which guarantees that $\boldsymbol{\varepsilon}'\neq\boldsymbol{\varepsilon}_+$, we have \begin{equation}\mathrm{val}_p(\#\tilde{\mathbb{T}}_{\Z_p}^{\boldsymbol{\varepsilon}'}/\tilde{J}_{\mathbb{Z}_p}^{\boldsymbol{\varepsilon}'}) = \mathrm{val}_p\Big(\prod_{j=1}^t(q_j+\varepsilon_j')\Big)=\val_p(\eta(N)).\end{equation}
Hence, since $\tilde{\mathbb{T}}_{\mathbb{Z}_p}/\tilde{J}_{\mathbb{Z}_p}\twoheadrightarrow\,\tilde{\mathbb{T}}_{\Z_p}^{\boldsymbol{\varepsilon}}/\tilde{J}_{\mathbb{Z}_p}^{\boldsymbol{\varepsilon}}$ for each $\boldsymbol{\varepsilon}\in\textit{\textbf{E}},$ we conclude
\begin{equation}\label{ohtap}\mathrm{val}_p(\#\tilde{\mathbb{T}}_{\mathbb{Z}_p}/\tilde{J}_{\mathbb{Z}_p})\geq\mathrm{val}_p(\eta(N)).\end{equation}
To obtain the desired lower bound 
\begin{equation}\label{ohtabound}\mathrm{val}_p(\#\mathbb{T}_{\mathbb{Z}_p}/J_{\mathbb{Z}_p})\geq\mathrm{val}_p(\eta(N)),\end{equation} we observe that the natural map $\mathbb{T}_{\mathbb{Z}_p}/J_{\mathbb{Z}_p}\to\, \tilde{\mathbb{T}}_{\mathbb{Z}_p}/\tilde{J}_{\mathbb{Z}_p}$ is in fact a surjection since the coset of $\mathbb{T}_{\mathbb{Z}_p}/J_{\mathbb{Z}_p}$ containing $1+p$ maps onto the coset of $\tilde{\mathbb{T}}_{\mathbb{Z}_p}/\tilde{J}_{\mathbb{Z}_p}$  containing $T_p$. We summarize these results in the following proposition:

\begin{proposition}Let $N$ be a squarefree integer and $p\geq 5$ a prime that does not divide $N$. If none of the prime divisors of $N$ are congruent to $1\pmod p$ and at least one prime divisor of $N$ is not congruent to $-1\pmod p,$ then there is an equality
\[\val_p(\eta(N))=\val_p(\#\mathbb{T}_{\Z_p}/J_{\Z_p}).\]
\end{proposition}
\begin{remark} This equality is important in the study of the modularity of residual Galois representations. Specifically, the lower bound in Eq. (\ref{ohtabound}) can be combined with the non-principality result in Theorem \ref{princondition} to give a statement analogous to \cite[Theorem 5.12]{berger2018modularity}, regarding the existence of many modular Galois representations, in the case of weight 2 cusp forms of trivial Nebentypus. Note that this context specifically requires the exclusion of the Hecke operator $T_p$ from the Hecke algebra.
\end{remark}

\begin{remark}\label{lowerbound} The results of Ohta used here do not require that each $q_j$ satisfies $q_j\not\equiv 1\pmod p$. Rather, as long as $N$ has a prime divisor that is not congruent to $-1\pmod p,$ Eq. (\ref{ohtap}) holds so that by Proposition \ref{congruenceresult}, the depth of Eisenstein congrunces mod $p$ is bounded from below by the $p$-adic valuation of $\eta(N),$ i.e.,
\[\frac{1}{e_N}(m_1+
\cdots+m_r)\geq\mathrm{val}_p(\eta(N)).\]

\end{remark}

\subsubsection{Counting Eisenstein congruences} Theorem \ref{princondition} gives a numerical criterion for $\Jz$ to be non-principal. While this allows us to use direct computations to find examples where the Eisenstein ideal is not locally principal, which we do in Section \ref{examples and application}, it would also be useful to find conditions on the level $N$ which suffice to show the associated Eisenstein ideal is not locally principal. By combining a counting argument with a result of Ribet--Yoo, we give one such condition below.

We first prove a lower bound on the number of newforms satisfying an Eisenstein congruence (away from $Np$) modulo $p$. Note that through the end of this section, we use the term Eisenstein congruence to mean a congruence away from $Np$.
\begin{theorem}\label{countingtheorem}
Let $N =\prod_{i=1}^t q_i$ be a squarefree integer and $p\geq 5$ be prime. Assume $\val_p(q_i+1)>0$ for $i = 1,\dots,s,$ where $1\leq s\leq t$.
There are at least $s\cdot 2^{t-2}$ newforms of level dividing $N$ that satisfy an Eisenstein congruence modulo $p$.
\end{theorem}
\begin{proof}
We require a result of Ribet--Yoo \cite[Theorem 1.3(3)/Theorem 2.2(2)]{yoo2014non}, which gives necessary and suffiencient conditions for the existence of Eisentein congruences. Indeed, the result of Ribet--Yoo is phrased in terms of Galois representations and admissible tuples of primes; we now restate it in terms of congruences:
\begin{proposition}[Ribet--Yoo, 2018]\label{Ribet-Yoo} Let $M = \prod_{j=1}^vr_j$ be a squarefree integer and $p\geq 5$ be prime. If $v$ is even and $r_v\equiv -1\pmod p,$ then there exists a newform $f$ of level $M$ such that $f$ satifies an Eisenstein congruence modulo $p$ and such that \begin{equation}\label{plusminuseigen} T_{r_j}f = \begin{cases}f & \text{ if }j=1,\dots v-1,\\-f& \text{ if }j=v,\end{cases} \end{equation} where $T_{r_j}$ denotes the usual Hecke operator. 

\end{proposition}

To obtain the lower bound in Theorem 1, we find, for each $q_i$ with $1\leq i\leq s,$ a set of $2^{t-2}$ newforms, each of which has level dividing $N$ and satisfies an Eisenstein congruences modulo $p$. We then show that these sets are disjoint.

Without loss of generality, consider $q_1\equiv -1\pmod p$. We apply Proposition \ref{Ribet-Yoo} with each divisor $M$ of $N$ such that $M$ is divisible by $q_1$ and $M$ is the product of an even number of prime divisors. Specifically, for each odd integer $n\leq t-1,$ there are $\binom{t-1}{n}$ choices for a divisor $M$ of $N$ that satisfies the required conditions, and so summing over choices of $n$ gives a total of \[\sum_{\substack{n\,\leq\, t-1\\n\,\,\mathrm{ odd}}} \textstyle{\binom{t-1}{n}} = 2^{t-2}\] 
choices for $M$. For each choice of $M$, we apply Proposition \ref{Ribet-Yoo} with $r_v = q_1\equiv -1 \pmod p$ to obtain a newform of level $M$ that satisfies an Eisenstein congruence modulo $p$. In particular, since each choice of $M$ is distinct, the multiplicity one theorem guarantees that these $2^{t-2}$ newforms will also be distinct. 

It remains to show that the sets of newforms associated to the $q_i$ (with $1\leq i\leq s$) are disjoint. Again, without loss of generality, suppose that $f_1$ (resp. $f_2$) is a newform associated to $q_1$ (resp. $q_2$). If the levels of the newforms $f_1$ and $f_2$ are not equal, then $f_1\neq f_2$ by multiplity one, as above. If the levels are equal, then $f_1\neq f_2$ since $T_{q_1} f_1 =-f_1$ but $T_{q_1} f_2 =f_2$ (by Eq. (\ref{plusminuseigen}) and our choices of $r_v$).
\end{proof}

Now, the conditions on the Hecke eigenvalues in Eq. (14) actually do more than distinguish the newforms obtained from Theorem \ref{countingtheorem}; they show that none of these $s\cdot 2^{t-2}$ newforms are Galois conjugates. In particular, as explained in Remark \ref{galoisremark}, this means that the $s\cdot 2^{t-2}$ newforms obtained from Theorem \ref{countingtheorem}, along with their Galois conjugates under the action of the appropriate decomposition group, contribute at least $s\cdot 2^{t-2}$ to the depth of Eisenstein congruences modulo $p,$ i.e., \[\frac{1}{e_N}(m_1+
\cdots+m_r)\geq s\cdot 2^{t-2}.\] Thus, combining this inequality with Theorem \ref{princondition} establishes Theorem \ref{nonprincthm}, which states that under the assumptions of Theorem \ref{countingtheorem}, if $s\cdot 2^{t-2}>\val_p(\eta(N)),$ then $J_\Z$ is not locally principal.

\begin{remark} Other results of Ribet--Yoo \cite{yoo2014non} and independent work of Martin \cite{martin2016jacquet} give more conditions (in the style of Proposition \ref{Ribet-Yoo}) for the existence of Eisenstein congruences mod $p$. One could use these conditions, within the higher congruences framework, to give additional statements similar to Theorem \ref{nonprincthm}.

\end{remark}

\begin{remark}\label{rmkwiththree} Since we have assumed $p\geq 5$ throughout this paper, we now briefly address the cases $p=2,3$. Indeed, when $p=2,3,$ the higher congruences framekwork established in Section \ref{framekwork}, including Proposition \ref{congruenceresult}, still holds. However, problems arise when we try to use this framework to determine whether the Eisenstein ideal is locally principal. Specifically, when $p=2,$ Proposition \ref{indexBK} is not applicable since every prime is congruent to $\pm 1\pmod 2$. When $p=3,$ the admissibility result of Ribet--Yoo recalled in Proposition \ref{Ribet-Yoo} is not applicable, and so Theorems \ref{nonprincthm} and \ref{countingtheorem} do not hold. Nonetheless, the criterion given in Theorem \ref{princondition} is still valid, and so when $p=3,$ we use direct computations (in Section \ref{examples and application}) to give examples where the associated Eisenstein ideal is not locally principal.

\end{remark}

\section{Applications and examples}\label{examples and application}

For squarefree level $N,$ Proposition \ref{congruenceresult} bounds
the depth of Eisenstein congruences modulo $p$ by the $p$-adic valuation of $\#\mathbb{T}_{\Z_p}/J_{\Z_p}$. In this section, we first express this depth of congruence as the multiplicity \[\frac{1}{e_N}\sum_{i=1}^rm_i=e(\Jz,\Tz).\] We then use MAGMA \cite{bosma1997magma} to give computational examples of our main results.

\subsection{Hilbert--Samuel multiplicities and elliptic modular forms} We apply the commutative algebra result stated in Corollary \ref{length} in the context of elliptic modular forms to obtain the following proposition:

\begin{proposition}\label{mainresult}
For $i=1,\dots,r,$ let $\varpi_N^{m_i}$ be the highest power of $\varpi_N$ such that the Hecke eigenvalues of $f_i$ are congruent to those of $E_{2,N}$ modulo $\varpi_N^{m_i}$ for Hecke operators $T_\ell$ for all primes $\ell\nmid Np$.   Then
\[\frac{1}{e_N}\sum_{i=1}^r m_i =e(\Jz,\Tz).\]
\end{proposition}
\begin{proof}
As in the proof of Proposition \ref{congruenceresult}, take $T = \To$ and $J=\Jo,$ where $\mathfrak{m}$ is the unique maximal ideal of $T_\Oo$ containing $J_\Oo$. Let $T_i=\Oo$ and $\varphi_i:T\to\,T_i$ be the map sending a Hecke operator to its eigenvalue corresponding to $f_i$. Also, let $\varphi_i(\Tz)=T_{{\Z_p},i},\,\varphi_i(\Jz)=J_{{\Z_p},i}$. We then have\[\begin{split}
\length_\Oo(T_i/J_i) &= \val_\varpi(\#T_i/J_i) \\&= \frac{1}{d_N}\cdot\val_p(\#T_i/J_i)\\&=\frac{[\Oo:\Z_p]}{d_N}\cdot\val_p(\#T_{{\Z_p},i}/J_{{\Z_p},i})\\&=e_N\cdot\length_{\Z_p}(T_{{\Z_p},i}/J_{{\Z_p},i}).
\end{split}\]
Since $J_i,$ and $J_{{\Z_p}_i}$ are principal for each $i,$ we may apply Corollary \ref{length} to obtain \[
e(\Jo,\To) = e_N\cdot e(\Jz,\Tz).\]
Thus, 
\[\sum_{i=1}^s\val_p\left(\#T_i/J_i\right) =d_N\cdot e(\Jo,\To)
= [\Oo:\Z_p]\cdot e(\Jz,\Tz),\]
and combining this with Eq. (\ref{TiJi}) yields the desired result.
\end{proof}

\subsection{Computational examples}

We compute Eisenstein congruences for a selection of squarefree levels. To keep these computations to a manageable size, we actually compute congruences away from $N$ rather than $Np$. However, since a congruence away from $N$ is necessarily a congruence away from $Np,$ these computations suffice to show that the Eisenstein ideal $J_{\mathbb{Z}_p}$ is not locally principal. In fact, since Theorem \ref{princondition} is applicable in the case of the Eisenstein ideal $\tilde{J}_{\Z_p},$ which includes the generator $T_p-(1+p),$ these computations also establish examples where the local Eisenstein ideal $\tilde{J}_{\mathcal{\Z}_p,\mathfrak{m}_{\Z_p}}$ is non-principal.

Recall from Section 1 that we want to compute congruences between the Hecke eigenvalues of weight 2 newforms  $f_1,\dots,f_r$ of level $N_{f_i}$ dividing $N$ and the weight 2 Eisenstein series $E_{2,N}$. Since these forms are normalized eigenforms for all Hecke operators $T_\ell$ with $\ell\nmid N$ prime, this is equivalent to computing congruences between Fourier coefficients, i.e., congruences of the type \begin{equation}\label{congruencer}a_\ell(f_i)\equiv a_\ell(E_{2,N})\pmod{{\lambda_i}^r},\end{equation}
for all primes $\ell\nmid N$. While the algorithm we use is discussed in more detail in Appendix \ref{algorithm}, we give a sample data entry and a brief explanation below.

\begin{center}\vspace{.2cm}
$N=78 = 2\times 3\times 13,$ $p=7$, $\mathrm{val}_p(\eta(N))=1$:
\[\begin{array}{|c|c|c|c|c|} \hline\mathrm{level} &\mathrm{depth} & 
\mathrm{ramindex} & \mathrm{resfield} & \mathrm{conjclass} \\\hline
26&1&1&7&2\\
39&1&1&7&2\\
\hline\end{array}\]
\end{center}
Each line of this table corresponds to a newform $f_i$ that represents its Galois orbit under $\mathrm{Gal}(\overline{\Q}/\Q)$. Column 1 gives the level $N_{f_i}$ of $f_i,$ and Column 5 gives the number of the Galois orbit of $f_i$ with respect to the internal MAGMA numbering. For each congruence, $\lambda_i$ is a prime ideal, above the prime $p\in\Z,$ in the ring of integers of the coefficient field $K_{f_i}$. Column 2 gives the exponent of each congruence, i.e., the value of $r$ in Eq. (\ref{congruencer}), and Columns 3 and 4 give the ramification index and the order of the residue field, respectively, of the ideal $\lambda_i$ at $p$. Note that to simplify calculations, we compute each congruence in the ring of integers of the individual coefficient field $K_{f_i}$. In particular, because ramification indices are multiplicative, we can easily translate this data into congruences modulo a uniformizer of the ring of integers in the composite coefficient field $E/\Q_p,$ as required by Proposition \ref{congruenceresult}.

\begin{remark} As discussed in Section \ref{depthsresult}, Theorems \ref{nonprincthm} and \ref{princondition} give numerical criteria that can be used to show that the Eisenstein ideal $J_{\Z}$ is not locally principal. By combining explicit computations with the multiplicity result in Proposition \ref{mainresult}, we might be able to give a more precise bound on the minimal number of generators that each (local) Eisenstein ideal $\Jz$ requires. 
\end{remark}

\subsubsection{Examples where the Eisenstein ideal is not locally principal}
We use direct computations and the numerical criterion in Theorem \ref{princondition} to give examples where the Eisenstein ideal is not locally principal. Note that in each of these examples, the integers $p$ and $N$ satisfy the assumptions of Proposition \ref{indexBK}, i.e., $p\nmid N$ and $N$ has no prime factors congruent to $1\pmod p$.  Because the depth of Eisenstein congruences modulo $p$ is strictly greater than $\mathrm{val}_p(\eta(N)),$ we conclude that the (local) Eisenstein ideal $\Jz$ is not principal.

\begin{center}

\vspace{.2cm}
$N=195 = 3\times 5\times 13,$ $p=7$, $\mathrm{val}_p(\eta(N))=1$:
\[\begin{array}{|c|c|c|c|c|} \hline\mathrm{level} &\mathrm{depth} & 
\mathrm{ramindex} & \mathrm{resfield} & \mathrm{conjclass} \\\hline
39&1&1&7&2\\
65&1&1&7&3\\
\hline\end{array}\]

$N=354 = 2\times 3\times 59,$ $p=5$, $\mathrm{val}_p(\eta(N))=1$:
\[\begin{array}{|c|c|c|c|c|} \hline\mathrm{level} &\mathrm{depth} & 
\mathrm{ramindex} & \mathrm{resfield} & \mathrm{conjclass} \\\hline
118&1&1&5&3\\
177&1&2&5&2\\
\hline\end{array}\]

\vspace{.2cm}
$N=618 = 2\times 3\times 103,$ $p=13$, $\mathrm{val}_p(\eta(N))=1$:
\[\begin{array}{|c|c|c|c|c|} \hline\mathrm{level} &\mathrm{depth} & 
\mathrm{ramindex} & \mathrm{resfield} & \mathrm{conjclass} \\\hline
206&1&1&13&4\\
309&1&1&13&4\\
\hline\end{array}\]

\vspace{.2cm}
$N=786 = 2\times 3\times 131,$ $p=11$, $\mathrm{val}_p(\eta(N))=1$:
\[\begin{array}{|c|c|c|c|c|} \hline\mathrm{level} &\mathrm{depth} & 
\mathrm{ramindex} & \mathrm{resfield} & \mathrm{conjclass} \\\hline
262&1&1&11&5\\
393&1&2&11&5\\
\hline\end{array}\]

\end{center}

Now, from Theorem \ref{countingtheorem}, we expect at least $2^{3-2}=2$ congruences in each of the above examples, and our computations verify this expectation. The following examples illustrate Theorems \ref{nonprincthm} and \ref{countingtheorem} in more complex situations, such as when $N$ has more than $3$ prime divisors or more than $1$ prime divisor that is congruent to $-1\pmod p$. In each case, $s\cdot 2^{t-2}>\val_p(\eta(N)),$ and so $\Jz$ is non-principal.

\begin{center}\vspace{.2cm}
$N=798=2\times 3\times 7\times 19,$ $p=5$, $\mathrm{val}_p(\eta(N))=1$:
\[\begin{array}{|c|c|c|c|c|} \hline\mathrm{level} &\mathrm{depth} & 
\mathrm{ramindex} & \mathrm{resfield} & \mathrm{conjclass} \\\hline
38&1&1&5&2\\
57&1&1&5&3\\
133&1&2&5&2\\
798&1&2&5&13\\
\hline\end{array}\]

\vspace{.2cm}
$N=1066= 2\times 13\times 41,$ $p=7$, $\mathrm{val}_p(\eta(N))=2$:
\[\begin{array}{|c|c|c|c|c|} \hline\mathrm{level} &\mathrm{depth} & 
\mathrm{ramindex} & \mathrm{resfield} & \mathrm{conjclass} \\\hline
26&1&1&7&2\\
82&1&1&7&2\\
533&1&2&7&3\\
533&1&2&7&5\\
1066&1&1&7&10\\
\hline\end{array}\]

\vspace{.2cm}
$N=1102= 2\times 19\times 29,$ $p=5$, $\mathrm{val}_p(\eta(N))=2$:
\[\begin{array}{|c|c|c|c|c|} \hline\mathrm{level} &\mathrm{depth} & 
\mathrm{ramindex} & \mathrm{resfield} & \mathrm{conjclass} \\\hline
38&1&1&5&2\\
58&1&1&5&2\\
551&1&1&5&7\\
551&1&1&5&8\\
1102&1&1&5&14\\
\hline\end{array}\]

\end{center}

\subsubsection{Examples where the Eisenstein ideal $\tilde{J}_{\Z_p,\mathfrak{m}_{\Z_p}}$ is principal} Using direct computations and Eq. (\ref{ohtap}), we can give examples of squarefree levels $N$ where the (local) Eisenstein ideal $\tilde{J}_{\Z_p,\mathfrak{m}_{\Z_p}}$ is principal. Our choices for $p$ and $N$ satisfy the conditions in Remark \ref{lowerbound}, i.e., $p\geq 5$ and $N$ has at least one prime divisor that is not congruent to $-1\pmod p$. In particular, we allow $N$ to have divisors which are congruent to $1\pmod p$. In each example, the depth of Eisenstein congruences is equal to $\val_p(\eta(N))$ so that by Proposition \ref{congruenceresult} and Eq. (\ref{ohtap}), $\tilde{J}_{\Z_p,\mathfrak{m}_{\Z_p}}$ is principal.

\begin{center}

\vspace{.2cm}
$N=145 = 5\times 29,$ $p=7$, $\mathrm{val}_p(\eta(N))=1$:
\[\begin{array}{|c|c|c|c|c|} \hline\mathrm{level} &\mathrm{depth} & 
\mathrm{ramindex} & \mathrm{resfield} & \mathrm{conjclass} \\\hline
29&1&1&7&1\\
\hline\end{array}\]

\vspace{.2cm}
$N=413 = 7\times 59,$ $p=5$, $\mathrm{val}_p(\eta(N))=1$:
\[\begin{array}{|c|c|c|c|c|} \hline\mathrm{level} &\mathrm{depth} & 
\mathrm{ramindex} & \mathrm{resfield} & \mathrm{conjclass} \\\hline
413&1&1&5&6\\
\hline\end{array}\]

\vspace{.2cm}
$N=515 = 5\times 103,$ $p=13$, $\mathrm{val}_p(\eta(N))=1$:
\[\begin{array}{|c|c|c|c|c|} \hline\mathrm{level} &\mathrm{depth} & 
\mathrm{ramindex} & \mathrm{resfield} & \mathrm{conjclass} \\\hline
515&1&1&13&4\\
\hline\end{array}\]

\vspace{.2cm}
$N=655=5\times 131,$ $p=11$, $\mathrm{val}_p(\eta(N))=1$:
\[\begin{array}{|c|c|c|c|c|} \hline\mathrm{level} &\mathrm{depth} & 
\mathrm{ramindex} & \mathrm{resfield} & \mathrm{conjclass} \\\hline
655&1&1&11&5\\
\hline\end{array}\]

\end{center}

\subsubsection{Examples with $p=3$} As discussed in Remark \ref{rmkwiththree}, Theorems \ref{nonprincthm} and \ref{countingtheorem} do not hold for $p=3$ because the admissibility result of Ribet--Yoo is invalid. However, we can still use the numerical criterion in Theorem \ref{princondition} and direct computations to give examples where the Eisenstein ideal is not locally principal. Note that when $p=3,$ we actually compare the depth of Eisenstein congruences to \[\val_p(B_2\cdot\eta(N)) = \mathrm{val}_p(\eta(N))-1,\] where $B_2=\frac{1}{6}$ denotes the second Bernoulli number. This correction factor appears in the general verion of the Berger--Klosin result \cite[\S 5.1, Proposition 5.6]{berger2018modularity}.

\begin{center}
\vspace{.2cm}
$N=110 = 2\times 5\times 11,$ $p=3$, $\mathrm{val}_p(\eta(N))=3$:
\[\begin{array}{|c|c|c|c|c|} \hline\mathrm{level} &\mathrm{depth} & 
\mathrm{ramindex} & \mathrm{resfield} & \mathrm{conjclass} \\\hline
110&1&1&3&1\\
110&1&1&3&2\\
110&1&2&3&4\\
\hline\end{array}\]

\vspace{.2cm}
$N=374 = 2\times 11\times 17,$ $p=3$, $\mathrm{val}_p(\eta(N))=4$:
\[\begin{array}{|c|c|c|c|c|} \hline\mathrm{level} &\mathrm{depth} & 
\mathrm{ramindex} & \mathrm{resfield} & \mathrm{conjclass} \\\hline
34&1&1&3&1\\
187&1&1&3&2\\
374&1&1&3&2\\
374&1&1&3&3\\
374&1&1&3&4\\
\hline\end{array}\]

\vspace{.2cm}
$N=935=5\times 11\times 17,$ $p=3$, $\mathrm{val}_p(\eta(N))=4$:
\[\begin{array}{|c|c|c|c|c|} \hline\mathrm{level} &\mathrm{depth} & 
\mathrm{ramindex} & \mathrm{resfield} & \mathrm{conjclass} \\\hline
85&1&2&3&3\\
187&1&1&3&2\\
935&1&1&3&2\\
935&1&1&3&8\\
935&1&1&3&9\\
\hline\end{array}\]

\end{center}

\appendix

\section{Algorithm for computations}\label{algorithm}

We give the algorithm implemented in MAGMA \cite{bosma1997magma} to compute the depth of Eisenstein congruences in Section \ref{examples and application}. This algorithm\footnote{The code used can be found at \url{https://people.maths.bris.ac.uk/~zx18363/research.html}.} has been adapted from \cite[$\S 4.2$]{naskrkecki2014higher}. Indeed, our main modification is to the Sturm bound in \cite[Theorem 2]{naskrkecki2014higher}:

\begin{lemma}\label{sturmbound}
Let $N$ be a positive integer, and let $f\in M_2(\Gamma_0(N))$ be a modular form with coefficients in $\mathcal{O}_K$ for some number field $K$. Let $\mathfrak{p}$ be a fixed prime lying over some rational prime $p,$ and suppose the Fourier coefficients of $f$ satisfy \[a_\ell(f)\equiv 0\mod{\mathfrak{p}^m}\]
for all primes $\ell\leq \mu'/6$ with $\ell\nmid N,$ where \[\mu' =[\mathrm{SL}_2(\Z):\Gamma_0(N')]\;\text{  for } N' =N\cdot\prod_{p|N}p. \]
Then $a_\ell(f)\equiv 0\mod{\mathfrak{p}^m}$ for all primes $\ell\nmid N$.
\end{lemma}
\begin{proof} Apply \cite[Lemma 4.6.5]{miyake2006modular} to obtain a modular form $f'\in M_2(\Gamma_0(N'))$ defined by \[f':=\sum_{\gcd{(n,N)}=1}a_n(f)\cdot q^n.\] Note that $N'=N\cdot\prod_{p|N}p$ as above. Since the Fourier coefficients of $f'$ are multiplicative for all $n$ such that $\gcd(n,N) =1$ and vanish at any $n$ such that $\gcd(n,N)\neq 1,$ the hypotheses of this lemma imply that \[a_n(f')\equiv 0\mod{\mathfrak{p}^m}\] for all $n\leq \mu'/6$. Hence, by the straightforward generalization of Sturm's theorem stated in \cite[Proposition 1]{chen2010congruences}, we have $f'\equiv 0\mod{\mathfrak{p}^m},$ and hence, \[a_\ell(f)\equiv 0\mod{\mathfrak{p}^m}\] for all primes $\ell\nmid N$.
\end{proof}

By Lemma \ref{sturmbound}, it is sufficient for our algorithm to check only for congruences between the Hecke eigenvalues of newforms $f_1,\dots,f_r$ and Eisenstein series $E_{2,N}$ for Hecke operators $T_\ell$ for primes $\ell\leq\mu'/6$ with $\ell\nmid N$. We therefore replace the Sturm bound in Naskr\k{e}cki's  algorithm with 
\[
B = \frac{1}{6} \cdot [\mathrm{SL}_2(\Z) : \Gamma_0(N')] = \frac{1}{6} \cdot  N \cdot \prod_{p|N}(p+1),
\] 
and check for congruences only at primes less than $B$.  
Since the utilization of orders in number fields in Naskr\k{e}cki's computations of congruences is unrelated to whether or not the level $N$ is prime, this adjusted Sturm bound allows us to generalize Naskr\k{e}cki's algorithm:
\newline\phantom{a}
\newline\noindent\textbf{Input:} A positive squarefree integer $N$.
For each non-prime divisor $M$ of $N$:
\begin{enumerate}[1.]
\item Compute Galois conjugacy classes of newforms in $S_2(\Gamma_0(M))$. Call the set $\textit{New}$.
\item Compute the Sturm bound \[B = \frac{1}{6}\cdot[\mathrm{SL}_2(\Z):\Gamma_0(N')] = \frac{1}{6}\cdot N\cdot\prod_{p|N}(p+1).\]
\item Compute the coefficients $a_\ell(E_{2,N})$ for primes $\ell\leq B$ with $\ell\nmid N$.
\item Calculate the set of primes $P=\left\{p\text{ prime}:p\,|\,\text{Numerator}\left(\eta(N)\right)\right\}$.
\item For each pair $(p,f)\in P\times\textit{New},$ compute $K_f,$ the coefficient field of $f$.
\item Find an algebraic integer $\theta$ such that $K_f = \Q(\theta)$.
\item Compute a $p$-maximal order $\mathcal{O}$ above $\Z[\theta]$.
\item Compute the set $\mathcal{S} = \{\lambda\in\mathrm{Spec}\mathcal{O}:\lambda\cap\Z = p\Z\}$.
\item For each $\lambda\in\mathcal{S},$ compute \[r_\lambda = \min_{\substack{\ell\,\,\,\mathrm{prime}\\\ell\leq B,\,\,\ell\nmid N}}\left(\mathrm{ord}_\lambda(a_\ell(f)-a_\ell(E_{2,N}))\right).\]
\end{enumerate}
\textbf{Output:} If $r_\lambda>0,$ then we have a congruence \[a_\ell(f)\equiv a_\ell(E_{2,N})\mod{(\lambda\mathcal{O}_f)^{r_\lambda}}\] for all primes $\ell\nmid N$.

\begin{remark}\label{galoisremark} Since this algorithm computes congruences modulo prime ideals in the ring of integers of a global field, we must reinterpret its output within the local framework used in Proposition \ref{congruenceresult}. More specifically, let $f_1,\dots,f_r$ be all newforms of level $M,$ and let $L/\Q$ contain all Fourier coefficients of the $f_i$'s. If $\mathfrak{p}\subseteq\mathcal{O}_L$ corresponds to our choice of embedding $\overline{\Q}_p\hookrightarrow \mathbb{C},$ then Proposition \ref{congruenceresult} requires us to check for congruences modulo $\mathfrak{p}$ for every $\mathrm{Gal}(L_\mathfrak{p}/\Q_p)$-orbit in the set of newforms $\{f_1,\dots,f_r\}$. Our algorithm accomplishes this by fixing one representative of each $\mathrm{Gal}(L/\Q)$-orbit in $\{f_1,\dots,f_r\}$ and checking for congruences modulo all prime ideals in $\mathcal{O}_L$ lying over $p$. Because the depth of Eisenstein congruences is scaled by the ramification index $e(\mathfrak{p}/p),$ each $\mathrm{Gal}(L_\mathfrak{p}/\Q_p)$-orbit of newforms will contribute to the depth of congruence at least a total equal to the residue degree $[\mathcal{O}_L/\mathfrak{p}:\Z_p/(p)]$. 
\end{remark}

\end{document}